\documentclass[12pt]{article}

\def\wt{\mbox{wt}}

\def\dd{\mbox{dd}}

\title{Deletion-correcting codes and dominant vectors }

\author{Emil Kolev\\ Institute of Mathematics and Informatics \\ Bulgarian Academy of Sciences\\
e-mail: emil@math.bas.bg}
\begin{document}

\maketitle

\abstract{In this paper we describe all pairs of binary vectors $({\bf u}, {\bf v})$ such that the set of vectors obtained by $t$ deletions in ${\bf v}$ is a subset of the set of vectors obtained by $t$ deletions in ${\bf u}$ for $t=1,2$. Such pairs play an important role for finding the value of $L_2(n,t)$, the maximum cardinality of binary $t$-deletion-correcting code of length $n$}



\section{Introduction}

When a binary message is transmitted through a noisy channel some of its symbols may change. The receiver needs reliable tools for recovering the message. This is done by adding some extra symbols (called check symbols) to the original message and the result is a longer message.  The set of all such messages is called an error-correcting code. One of the main goals of coding theory is finding codes with good error-correcting capabilities.

Another possible distortion of the message is the lost of some of its symbols or insertion of some extra symbols. In this case the receiver gets shorter or longer message and he does not know which of the symbols were lost or inserted. Deletion-correcting codes and insertion-correcting codes are designed to correct such deletions or insertions. Levenstein has shown \cite{Lev65} that deletion-correcting codes and insertion-correcting codes are essentially the same objects. In this paper we consider only deletion-correcting codes.
A code is called $t$-deletion-correcting if it corrects any $t$ deletions. For more information and useful results the reader is referred to \cite{HF02}, \cite{LH06}, \cite{Lev65}, \cite{Lev65a},  \cite{Lev92}, \cite{Lev01}, \cite{Slo00}, \cite{SF03}, \cite{Tol02}.

\smallskip

{\bf Example 1.} Consider the binary code ${\cal C}=\{00000,11111,00011,11000,10101, 01110\}$. For a given codeword we may delete any of its five symbols. As a result we obtain a set of vectors of length 4. Direct verification shows that all six sets obtained from the six codewords are disjoined. Therefore ${\cal C}$ is 1-deletion-correcting code.

\smallskip

{\bf Definition 1.} The {\it Levenstein distance} $d_L({\bf x},{\bf y})$ of two binary vectors is defined as the minimum number of deletions and insertions needed to transform ${\bf x}$ into ${\bf y}$.

For example, $d_L(0100,110101)=4$. Note that in the above definition the vectors ${\bf x}$ and ${\bf y}$ do not need to be of one and the same length.

{\bf Definition 2.} {\it Deletion distance} $\dd({\bf u},{\bf v})$ between two vectors ${\bf u}$ and ${\bf v}$ of equal length is defined as
one-half of the smallest number of deletions and insertions needed to change ${\bf u}$ to ${\bf v}$, \cite{Slo00}.

For example, $\dd(00000,11111)=5$ whereas $\dd(00011,10101)=2$.
It is clear that for vectors ${\bf u}$ and ${\bf v}$ of equal length we have
$$
\dd({\bf u},{\bf v})=\frac{1}{2}d_L({\bf u},{\bf v}).
$$
For a given code ${\cal C}$ the deletion distance $\dd({\cal C})$ is defined as
$$
\dd({\cal C})=\min\{\dd({\bf u},{\bf v})\ | \ {\bf u},{\bf v}\in {\cal C}\}.
$$
For any two distinct codewords ${\bf u}$ and ${\bf v}$ from $t$-deletion-correcting code ${\cal C}$ of length $n$ we have $\dd({\bf u},{\bf v})>t$ (or, equivalently $d_L({\bf u},{\bf v})>2t$).

Denote by $L_2(n,t)$ the maximum cardinality of a binary $t$-deletion-correcting code ${\cal C}$ of length $n$.
A binary $t$-deletion-correcting code ${\cal C}$ of length $n$ and cardinality $L_2(n,t)$ is called optimal.

For a binary vector ${\bf u}$ of length $n$ denote by $D_t({\bf u})$ the set of all vectors of length $n-t$ obtained from ${\bf u}$ by deleting $t$ entries in ${\bf u}$. In other words, $D_t({\bf u})$ contains all subsequences of ${\bf u}$ of length $n-t$.

The size of $D_t({\bf u})$ depends on ${\bf u}$. The minimal size of $D_t({\bf u})$ equals 1 and is achieved only for ${\bf u}=p^n$ for $p\in\{0,1\}$. The problem of finding the maximal size of $D_t({\bf u})$ is discussed in \cite{Calabi}, \cite{Lev2001}.

A code ${\cal C}$ is $t$-deletion-correcting code if the sets $D_t({\bf u})$ for $u\in {\cal C}$ are disjoint. Further, if the sets $D_t({\bf u})$ for $u\in {\cal C}$ partition the set $F_q^{n-t}$ then the code is called perfect.

As in the case of error-correcting codes the two main research problems for deletion-correcting codes are:

1. For given $n$ and $t$ find $L_2(n,t)$, the maximum cardinality of a binary $t$-deletion-correcting code of length $n$.

2. When $L_2(n,t)$ is known, find all distinct (in some sense) optimal codes.

\smallskip

In general, finding the value of $L_2(n,t)$ is an open problem in coding theory.
The efforts are concentrated on specific values of $n$ and $t$.
Tables with known values of $L_2(n,t)$ for different $n$ and $t$ can be found in \cite{LH06} and \cite{Kalin}.

\section{Preliminaries}

Any permutation of coordinates of given code ${\cal C}$ does not alter its error-correcting capabilities.
On the contrary, for deletion-correcting codes a permutation of coordinates, in general, does not result in a code with the same deletion-correcting properties. Nevertheless, there are two simple observations that describe when two deletion-correcting codes are essentially the same and allow to adopt different notion for equivalence. First, we may read the codewords backwards and second, we may change 0 and 1. This leads to the following

\smallskip

{\bf Definition 3.} Two deletion-correcting codes ${\cal C}_1$ and ${\cal C}_2$ are equivalent if
one of the following is true:

1. $(u_1,u_2,\dots ,u_n) \in {\cal C}_1$ if and only if  $(\overline{u_1},\overline{u_2},\dots ,\overline{u_n}) \in  {\cal C}_2$;

2. $(u_1,u_2,\dots ,u_n) \in {\cal C}_1$ if and only if  $(u_n,u_{n-1},\dots ,u_1) \in  {\cal C}_2$;

3. $(u_1,u_2,\dots ,u_n) \in {\cal C}_1$ if and only if  $(\overline{u_n},\overline{u_{n-1}},\dots ,\overline{u_1}) \in  {\cal C}_2$.

Here, for $x\in\{0,1\}$ the element $\overline{x}\in\{0,1\}$ is such that $\{x,\overline{x}\}=\{0,1\}$.

In finding the exact value of $L_2(n,t)$ usually at some stage an exhaustive computer search is performed.
As in any computer search a good pruning technique is required. It turns out that when choosing the codewords of optimal deletion-correcting code some of the vectors may be left out.


{\bf Definition 4.} We say that a vector ${\bf u}$ is $t$-dominant if there exists a vector ${\bf v}$ (alternatively, ${\bf v}$ is subordinate of ${\bf u}$) such that ${\bf u}\neq {\bf v}$ and $D_t({\bf v})\subseteq D_t({\bf u})$.

It is clear that if ${\bf u}$ is $t$-dominant over ${\bf v}$ then for any $s>t$ the vector ${\bf u}$ is $s$-dominant over the vector ${\bf v}$.
If a codeword ${\bf u}$ is $t$-dominant over the vector ${\bf v}$ then
$$
{\cal C}\setminus \{{\bf u}\} \cup \{{\bf v}\}
$$
is also $t$-deletion-correcting code. In other words a dominant codeword may be replaced by its subordinate vector. Hence, in computer search we may exclude all dominant vectors from consideration. Therefore it is important to know all pairs of vectors $({\bf u},{\bf v})$ such that $D_t({\bf v})\subseteq D_t({\bf u})$.

Furthermore, we may assume that an optimal code ${\cal C}$ includes the vectors $0^n$ and $1^n$ as codewords. Indeed, for $p\in\{0,1\}$:
\begin{itemize}
\item if $p^{n-t}\in D_t({\bf u})$ for a codeword ${\bf u}$ then, as above, replace ${\bf u}$ by $p^n$ and
\item if $p^{n-t}\not\in D_t({\bf u})$ for any codeword ${\bf u}$ then ${\cal C}\cup \{p^n\}$ is $t$-deletion-correcting code, i.e. ${\cal C}$ is not optimal.
\end{itemize}

A code ${\cal C}$ is called basic if it does not contain dominant vectors.
In the lights of the last two definitions the main problems for deletion-correcting codes become:

1. For certain $n$ and $t$ find $L_2(n,t)$;

2. Find all inequivalent basic optimal codes.

\section{Results}

As explained in the previous section knowing the pairs of $t$ dominant vectors plays an important role in finding $L_2(n,t)$.
In what follows we describe all pairs of binary vectors $({\bf u},{\bf v})$ such that ${\bf u}$ is $t$-dominant over ${\bf v}$ for $t=1$ and $t=2$.

For the two trivial cases ${\bf v}=0^n$, ${\bf v}=1^n$ and for any $t$ we have:
\begin{itemize}
\item if ${\bf v}=0^n$ then ${\bf u}$ is $t$-dominant over ${\bf v}$ if and only if ${\bf u}\neq {\bf v}$ and $\wt({\bf u})\leq t$;
\item if ${\bf v}=1^n$ then ${\bf u}$ is $t$-dominant over ${\bf v}$ if and only if ${\bf u}\neq {\bf v}$ and $\wt({\bf u})\geq n-t$.
\end{itemize}

In what follows the vector ${\bf u}=(u_1,u_2,\dots ,u_n)$ is $t$-dominant over ${\bf v}=(v_1,v_2,\dots ,v_n)$ and $\{p,q\}=\{0,1\}$.
We begin with a useful observation.

\smallskip

{\bf Proposition 1.} Let $n\geq 2$ be positive integer. Consider two vectors ${\bf x}$ and ${\bf y}$ of lengths $n$ and $n-1$ respectively. If any single deletion changes ${\bf x}$ to ${\bf y}$ then all entries in ${\bf x}$ and ${\bf y}$ are equal.

\smallskip

{\bf Proof.} Let ${\bf x}=(x_1,x_2,\dots ,x_n)$ and ${\bf y}=(y_1,y_2,\dots ,y_{n-1})$ and choose a positive integer $k$ such that $1\leq k\leq n-1$.
By deleting $x_{k}$ we have that $x_{k+1}=y_k$ and by deleting $x_{k+1}$ we infer that $x_k=y_{k}$. Therefore $x_k=x_{k+1}=y_k$ for any $k=1,2,\dots ,n-1$ which implies that $x_1=x_2=\cdots =x_n=y_1=y_2=\cdots =y_{n-1}$. $\diamond$

{\bf Remark.} The above proposition is true also for vectors ${\bf x}$ and ${\bf y}$ of lengths $n\geq 3$ and $n-2$, respectively, when the result of any two deletions in ${\bf x}$ is ${\bf y}$. The proof is straightforward.

\smallskip

First, we describe all $1$-dominant vectors.

{\bf Proposition 2.} Let ${\bf u}$ be $1$-dominant over ${\bf v}$ and ${\bf v}\neq 0^n,1^n$. Then  ${\bf u}=p^{m-1}qpq^{n-m-1}$ and ${\bf v}=p^mq^{n-m}$ for some positive integer $m$.

\smallskip

{\bf Proof.} Let ${\bf u}=(u_1,u_2,\dots ,u_n)$ and ${\bf v}=(v_1,v_2,\dots ,v_n)$. For $n=2$ the result is trivial, so let $n\geq 3$.
Assume first that $u_1\neq v_1$ and let $v_1=p$, $u_1=q$. Any deletion of $v_i$ for $i\geq 2$ results in a vector ${\bf w}\in D_1({\bf u})$ with first coordinate $v_1\neq u_1$. This is possible only if ${\bf w}$ is obtained from ${\bf u}$ by deleting its first coordinate and $u_2=v_1=p$. Proposition 1 applied for ${\bf x}=(v_2,\dots ,v_n)$ and ${\bf y}=(u_3,\dots ,u_n)$ implies that
$v_2=v_3=\cdots =v_n=u_3=\cdots =u_n$. Since ${\bf v}\neq p^n$ we infer that
$$
 {\bf u}=qpq^{n-2} \mbox{ and } {\bf v}=pq^{n-1}. \leqno(1)
$$
It is easy to check that ${\bf u}$ is $1$-dominant over ${\bf v}$. In this case $m=1$.

Assume ${\bf u}=(p,\dots, p, u_{k+1}, \dots ,u_n)$ and ${\bf v}=(p,\dots, p, v_{k+1}, \dots ,v_n)$ where $k\geq 1$ and $u_{k+1}\neq v_{k+1}$.

If $v_{k+1}=q$ then $u_{k+1}=p$. By deleting the first coordinate in ${\bf v}$ we obtain a vector ${\bf w}$ with $k$-th coordinate equals to $q$. Note that all vectors from $D_1({\bf u})$ have their first $k$ entries equal to $p$. Therefore ${\bf w}\not\in D_1({\bf u})$, a contradiction.

Hence, $v_{k+1}=p$ and $u_{k+1}=q$. Since ${\bf v}\neq p^n$ we have that $n\geq k+2$.
By deleting $v_{i}$ for arbitrary $i\geq k+2$ we obtain a vector ${\bf w}$ with first $k+1$ entries equal to $p$.
The only way to obtain such a vector by 1 deletion in ${\bf u}$ is to have $u_{k+2}=p$ and to delete $u_{k+1}=q$. If $n=k+2$ then ${\bf u}=p^{n-2}qp$, ${\bf v}=p^{n-1}q$, and this pair is equivalent to the pair described in (1). If $n\geq k+3$ then  Proposition 1 applied for ${\bf x}=(v_{k+2},\dots ,v_n)$ and
${\bf y}=(u_{k+3},\dots ,u_n)$ implies that $v_{k+2}=\cdots =v_n=u_{k+3}=\cdots =u_n$. Since ${\bf v}\neq p^n$ we conclude that ${\bf v}=p^{k+1}q^{n-k-1}$ and ${\bf u}=p^{k}qpq^{n-k-2}$.  In this case $m=k+1$. $\diamond$

\smallskip

In Table 1 we present all pairs ${\bf u}$ and ${\bf v}$ such that ${\bf u}$ is $1$-dominant over ${\bf v}$.

\medskip

\begin{center}
$\begin{array}{|c|c|c|}
 \hline
&  {\bf u}  &{\bf v}\\\hline
1. &\wt({\bf u})=1& 0^n\\\hline
2. &\wt({\bf u})=n-1& 1^n\\\hline
3. &p^{m-1}qpq^{n-m-1} & p^mq^{n-m}\\\hline
\end{array}$

\bigskip

Table 1.
\end{center}

\bigskip

We proceed now with the case $t=2$. Since the case ${\bf v}=p^n$ is
clear in what follows we assume that ${\bf v}\neq p^n$.

For $n=3$ up to equivalence we have: ${\bf v}=ppq$ and ${\bf u}\neq p^3,q^3,{\bf v}$ or
${\bf v}=pqp$ and ${\bf u}\neq p^3,q^3,{\bf v}$.

For $n=4$ we have that up to equivalence there exist 5 choices for
${\bf v}$, namely: $pppq$, $ppqp$, $ppqq$, $pqpq$ and $pqqp$. For
any of these instances it is easy to enumerate all vectors ${\bf u}$
that are $2$-dominant over ${\bf v}$.

\smallskip

Let $n\geq 5$ be positive integer and ${\bf u}=(u_1,u_2,\dots ,u_n)$
be $2$-dominant over ${\bf v}=(v_1,v_2,\dots ,v_n)$. Denote
$k=\min\{i | u_i\neq v_i\}$ and $s=\max\{i | u_i\neq v_i\}$ where
$k\leq s$. We split the proof in several cases depending on $k$ and
$s$.

\medskip

{\bf Case A.} $k=s$, i.e. $d({\bf u},{\bf v})=1$;

\medskip

{\bf Case B.} $k=1$ and $s=n$, i.e. $u_1\neq v_1$ and $u_n\neq v_n$;

\medskip

{\bf Case C.} $k\neq s$ and $u_1=v_1$.

\medskip


We settle the above cases in the next three propositions.

\smallskip

{\bf Proposition 3.} If ${\bf u}$ is $2$-dominant over ${\bf v}$ and
$d({\bf u},{\bf v})=1$ then up to equivalence ${\bf
u}=p^mqp^{n-m-2}q$ and  ${\bf v}=p^{n-1}q$ or ${\bf
u}=p^mqp^{n-m-3}qp$ and ${\bf v}=p^{n-2}qp$ for some integer $m\geq
0$.

{\bf Proof.} Since $d({\bf u},{\bf v})=1$ we have that there exists
positive integer $k$ such that $u_i=v_i$ for $i\neq k$ and $u_k=q$,
$v_k=p$. The number of elements $q$ in ${\bf u}$ is one more than
the corresponding entries in ${\bf v}$. Therefore if there exist two
or more entries $q$ in ${\bf v}$ then the vector ${\bf w}$ obtained
by deleting two elements $q$ in ${\bf v}$ has at least three
elements $q$ less than ${\bf u}$. Therefore ${\bf w}\not\in D_2({\bf
u})$. Since ${\bf v}\neq p^n$ we infer that ${\bf v}=p^bqp^{n-b-1}$
for some integer $b$ for which $0\leq b\leq n-1$. Up to equivalence
${\bf u}=p^{k-1}qp^{b-k}qp^{n-b-1}$. If $n-b-1\geq 2$ then the
deletion of the last two symbols from ${\bf v}$ gives a vector not
in $D_2({\bf u})$. Thus, $n-b-1=0$ or 1 and we obtain ${\bf
u}=p^mqp^{n-m-2}q$ and ${\bf v}=p^{n-1}q$ or ${\bf
u}=p^mqp^{n-m-3}qp$ and ${\bf v}=p^{n-2}qp$. It is easy to check
that in both cases ${\bf u}$ is $2$-dominant over ${\bf v}$.
$\diamond$


{\bf Proposition 4.} Let ${\bf u}$ be $2$-dominant over ${\bf v}$
and $u_1\neq v_1$, $u_n\neq v_n$. Then up to equivalence ${\bf
u}=qpq^{n-4}pq$ and ${\bf v}=pq^{n-2}p$ or ${\bf u}=qp^{n-3}qp$ and
${\bf v}=p^{n-1}q$.

{\bf Proof.} Without loss of generality assume $v_1=p$ and $u_1=q$.

1. Let $v_n=p$ and  $u_n=q$. The deletion of any two elements from
$v_2,v_3,\dots ,v_{n-1}$ gives a vector from $D_2({\bf u})$ with
first coordinate $v_1=p$ and last coordinate $v_n=p$. Such a vector
can be obtained from ${\bf u}$ only if we delete $u_1=q$ and
$u_n=q$. Therefore $u_2=u_{n-1}=p$ and any two deletions from
$(v_2,v_3,\dots ,v_{n-1})$ imply $(u_3,u_4,\dots ,u_{n-2})$. It
follows from the remark after Proposition 2 that $v_2=v_3=\cdots
=v_{n-1}=u_3=u_4=\cdots =u_{n-2}$. Since ${\bf v}\neq p^n$ we have
that ${\bf u}=qpq^{n-4}pq$ and ${\bf v}=pq^{n-2}p$. Direct
verification shows that indeed ${\bf u}$ is $2$-dominant over ${\bf
v}$.

2. Let $v_n=q$ and  $u_n=p$. As in the previous case we infer that
$u_2=p$, $u_{n-1}=q$ and $v_2=v_3=\cdots =v_{n-1}=u_3=\cdots
=u_{n-2}$. Up to equivalence ${\bf u}=qp^{n-3}qp$, ${\bf
v}=p^{n-1}q$ it is easy to see that ${\bf u}$ is $2$-dominant over
${\bf v}$.  $\diamond$


{\bf Proposition 5.} Let ${\bf u}$ be $2$-dominant over ${\bf v}$,
$u_1=v_1$ and $k\neq s$ where $k=\min\{i | u_i\neq v_i\}$ and
$s=\max\{i | u_i\neq v_i\}$. Then up to equivalence all such vectors
${\bf u}$ and ${\bf v}$ are presented in the following table.

$$
\begin{array}{|c|c|c|}\hline
&{\bf u}& {\bf v} \\\hline 1. &
pqp^{m-1}qpq^{n-m-3}&pqp^mq^{n-m-2}\\\hline 2. &
pq^mpqp^{n-m-3}&pq^{m+1}p^{n-m-2}\\\hline 3. &
p^2qp^{n-3}&pqp^{n-2}\\\hline 4. &
p^{m-2}qppqp^{n-m-2}&p^mqp^{n-m-1}\\\hline 5. &
p^{m-1}qpqp^{n-m-2}&p^mqp^{n-m-1}\\\hline 6.&p^{n-4}qppq& p^{n-2}qp
\\\hline
7.& p^{n-3}qpq  &p^{n-2}qp \\\hline 8.& p^{m}qp^{n-m-3}qp &
p^{n-1}q\\\hline 9. &p^{m-2}qpqpq^{n-m-2}& p^mq^{n-m}\\\hline 10.
&p^{m-2}qppqq^{n-m-2}& p^mq^{n-m}\\\hline 
11. &p^{n-2}qp&
p^{n-1}q\\\hline 12. &p^{n-3}qp^2& p^{n-1}q\\\hline 
13.& p^{m-1}qpqpq^{n-m-3}& p^mqpq^{n-m-2}\\\hline 

14.& p^{n-3}q^2 p & p^{n-2}q^2 \\\hline
15. &p^{n-4}qppq & p^{n-3}q^2p\\\hline 
16.& p^{n-4}qpqp & p^{n-3}qpq \\\hline
 17.
&p^{m-1}q^{n-m-1}pq&p^{m}q^{n-m-1}p
\\\hline
18. &p^{m-1}q^{n-m}p&p^{m}q^{n-m}
\\\hline
\end{array}
$$

{\bf Proof.} Without lost of generality assume $v_1=u_1=p$. Note
that since $k\neq s$ we have $d({\bf u},{\bf v})>1$.

1. If $v_2=q$ and $u_2=q$ then the deletion of $v_1$ and an
arbitrary $v_i$ for $i\geq 3$ implies the deletion of $u_1$ in ${\bf
u}$. Thus ${\bf u}$ without its first $2$ entries is $1$-dominant
over ${\bf v}$ without its first $2$ entries. Hence, ${\bf u}=pq{\bf
u_1}$ and ${\bf v}=pq{\bf v_1}$ where ${\bf u_1}$ is $1$-dominant
over ${\bf v_1}$ and $d({\bf u_1},{\bf v_1})>1$. Therefore the pair
$({\bf u_1},{\bf v_1})$ is equivalent to one of the pairs from Table
1 and $d({\bf u_1},{\bf v_1})>1$. Only the third entry in Table 1
satisfies $d({\bf u_1},{\bf v_1})>1$. Hence, we obtain the following
pairs: ${\bf u}=pqp^{m-1}qpq^{n-m-3}$ and ${\bf v}=pqp^mq^{n-m-2}$;
${\bf u}=pq^mpqp^{n-m-3}$ and ${\bf v}=pq^{m+1}p^{n-m-2}$. In both
cases we have that ${\bf u}$ is $2$-dominant over ${\bf v}$.

2. If $v_2=q$ and $u_2=p$ then the deletion of $v_1$ and arbitrary
$v_i$ for $i\geq 3$ implies a vector ${\bf w}\in D_2({\bf u})$ with
first coordinate $v_2=q$. To obtain ${\bf w}$ from ${\bf u}$ by two
deletions we should have $u_3=q$ and we have to delete $u_1$ and
$u_2$. We infer that a single deletion in $(v_3,v_4,\dots ,v_n)$
gives $(u_4,u_5,\dots ,u_n)$. Proposition 1 implies that
$v_3=\cdots=v_n=u_4=\cdots =u_n$. Thus, ${\bf u}=ppqp^{n-3}$ and
${\bf v}=pqp^{n-2}$ or ${\bf u}=ppq^{n-2}$ and ${\bf v}=pq^{n-1}$.
For both pairs ${\bf u}$ is 2-dominant over ${\bf v}$ but only for
the first pair we have $d({\bf u},{\bf v})>1$.

Let $v_2=p$ and assume $p=v_1=v_2=\cdots =v_m\neq v_{m+1}=q$
for some $m\geq 2$. If $u_1=u_2=\cdots =u_m=p$ the deletion of $v_1$
and $v_2$ implies the deletion of two of the first $m$ elements in
${\bf u}$. Thus, ${\bf u}={\bf v}$, a contradiction.

We conclude that $k\leq m$ and then $u_1=\cdots =u_{k-1}=p$ and $u_k=q$.

For $p\in \{0,1\}$ and a vector ${\bf w}$ denote by $n_p({\bf w})$
the number of entries $p$ in the vector ${\bf w}$.

\begin{itemize}
\item If $n_q({\bf v})>n_q({\bf u})$ then delete two elements $p$ from ${\bf v}$
and let ${\bf w}$ be the resulting vector.
Since $n_q({\bf w})>n_q({\bf u})$ we infer that ${\bf w}\not\in D_2({\bf u})$.
\item If $n_q({\bf v})<n_q({\bf u})$ then $n_p({\bf v})>n_p({\bf u})$. If
$n_q({\bf v})\geq 2$, i.e. there exist at least two entries $q$ in
${\bf v}$ then we delete two elements $q$ from ${\bf v}$ and obtain
a contradiction as above. Therefore ${\bf v}=p^{m}qp^{n-m-1}$ and
$n_q({\bf u})=2$ (if $n_q({\bf u})\geq 3$ the deletion of a symbol
$p$ and the symbol $q$ in ${\bf v}$ gives a contradiction). If
$u_{m+1}=q$ then $k=s$, a contradiction. Assume first that
$n-m-1\geq 2$. If $u_i=q$ for $i\neq k$ and $i\leq m$ then, as
above, the deletion of $v_{n-1}=p$ and $v_n=p$ gives a
contradiction. Thus, for some $i>m+1$ we have $u_i=q$. It is easy to
see that up to equivalence there exist two choices for ${\bf u}$:
${\bf u}=p^{m-2}qppqp^{n-m-2}$ and ${\bf u}=p^{m-1}qpqp^{n-m-2}$. If
$n-m-1=1$ then ${\bf v}=p^{n-2}qp$, ${\bf u}=p^{n-4}qppq$ or ${\bf
u}=p^{n-3}qpq$ and if $n-m-1=0$ then ${\bf v}=p^{n-1}q$ and ${\bf
u}=p^{m}qp^{n-m-3}qp$.
\item Let $n_q({\bf v})=n_q({\bf u})$.
Note that in this case the deleted symbols from ${\bf v}$ are
identical to the deleted symbols from ${\bf u}$. If $m\geq k+2$ the
deletion of the first two entries in ${\bf v}$ gives a
contradiction. If $m=k+1$ then if there exists $i>m+1$ such that
$v_i=p$ then the deletion of $v_1$ and $v_i$ implies a
contradiction. Thus, ${\bf v}=p^mq^{n-m}$ and then for $n-m\geq 2$ we
have ${\bf u}=p^{m-2}qpqpq^{n-m-2}$ or ${\bf
u}=p^{m-2}qppqq^{n-m-2}$. For $n-m=1$ we have ${\bf v}=p^{n-1}q$ and
${\bf u}=p^{n-2}qp$ or ${\bf u}=p^{n-3}qp^2$.

Let $m=k$. If there exist at least two entries $p$ in
$(v_{m+2},\dots , v_n)$ then the deletion of these two elements
gives a contradiction.

If $u_n=v_n=p$ then we may show as above that $v_{n-1}=p$ and we
have at least two entries $p$ in $(v_{m+2},\dots , v_n)$, a
contradiction.

If $u_n=v_n=q$ then the same observations as above but starting from
right imply that ${\bf u}=p^{m-1}q{\bf w}pq^b$ and  ${\bf
v}=p^{m}q{\bf h}pq^{b+1}$. If ${\bf h}$ is not empty then the
deletion of any two elements from $q{\bf h}p$ implies the deletion
of $u_{m}$ and $u_{n-b}$. Proposition 1 implies that all entries in
$q{\bf h}p$ are equal which is not true. Therefore ${\bf h}$ is
empty and then ${\bf v}=p^{m}qpq^{n-m-2}$ and  ${\bf
u}=p^{m-1}qu_{m+1}u_{m+2}pq^{n-m-3}$. Since $n_q({\bf v})=n_q({\bf u})$ we
have that $\{u_{m+1},u_{m+2}\}=\{p,q\}$. Only one of the two cases gives
$2$-dominant vectors, namely ${\bf u}=p^{m-1}qpqpq^{n-m-3}$ and
${\bf v}=p^{m}qpq^{n-m-2}$.

It remains to consider the case $u_n\neq v_n$. If $n\leq m+3$ then an easy enumeration gives:
\begin{itemize}
\item ${\bf u}=p^{n-2}qp$ and ${\bf v}=p^{n-1}q$ for $n=m+1$;
\item ${\bf u}=p^{n-3}q^2 p$ and ${\bf v}=p^{n-2}q^2$ for $n=m+2$; 
\item ${\bf u}=p^{n-4}qppq$ and ${\bf v}=p^{n-3}q^2p$; ${\bf u}=p^{n-4}qpqp$ and ${\bf v}=p^{n-3}qpq$ for $n=m+3$.
\end{itemize}
If $n\geq m+4$ then  any two deletions in $(v_{m+1},\dots ,v_{n-1})$
imply  $u_{n-1}=v_n$ and the deletion of $u_{m}$ and $u_{n}$.
Proposition 1 implies that $q=v_{m+1}=\cdots =v_{n-1}=u_{m+2}=\cdots
=u_{n-1}$, thus ${\bf v}=p^{m}q^{n-m-1}v_n$ and ${\bf
u}=p^{m-1}q^{n-m-1}v_nu_n$.

Both choices of $u_n\neq v_n$ give $2$-dominant pair. $\diamond$

\end{itemize}

\end{document}